\documentclass[12pt]{article}

\usepackage{amssymb,amsmath,latexsym}
\usepackage[utf8]{inputenc}
\usepackage{textcomp}
\usepackage[T1]{fontenc}
\usepackage{authblk}
\usepackage{amssymb,amscd,amsmath,amsfonts}
\usepackage{latexsym}
\usepackage{amsthm}
\usepackage{enumerate}
\usepackage{verbatim}
\usepackage{graphicx}
\usepackage[dvips]{epsfig}
\usepackage{subcaption}
\usepackage{amsbsy}
\usepackage{mathtools}

\usepackage{inputenc}

\allowdisplaybreaks

\newtheorem{thm}{Theorem}[section]
\newtheorem{lem}[thm]{Lemma}

\newtheorem{dfn}[thm]{Definition}
\newtheorem{exa}[thm]{Example}

\newtheorem{rk}[thm]{Remark}

\numberwithin{equation}{section}

\newcommand{\C}{\mathbb C}
\newcommand{\R}{\mathbb R}

\newcommand{\N}{\mathbb N}

\newcommand{\e}{\varepsilon}

\begin{document}

\title{Hyers-Ulam stability of the first order difference equation with average growth rate}

\author
{Young Woo Nam 
}

\affil[]{\small  
      College of Engineering,
      Konkuk University, 05029 
      Seoul, Korea \\ \texttt{namyoungwoo@konkuk.ac.kr}
      }
\date{}

\maketitle

\begin{abstract}
The first order difference equation induced by the sequence of maps on $ \mathbb{C} $ has Hyers-Ulam stability where the limit of the geometric average of growth rate is convergent and not equal to one. 
We show no Hyers-Ulam stability where the average growth rate is (pre)periodic even though each periodic growth rate is strictly less than one. Examples of difference equation generated by time dependent maps which contains contracting maps and expanding maps are given. 
\end{abstract}
\footnote{ 2010 {\em Mathematical Subject Classification}: Primary 39A20; Secondary 39A45. \\
  {\em Key words and phrases}:
Hyers-Ulam stability; root test; difference equation; average growth rate.}

\section{Introduction and preliminaries}

Hyers-Ulam stability has been an approximation theory of functional equations, differential equations for decades. See \cite{hyers}, \cite{HJL}. 
Hyers-Ulam stability of difference equation has been studied in some researchers, \cite{AO1, AOR, BBP, BCL, nam1, Po}. Many of these topics are concentrated on the difference equation generated by a single map. For instance, contracting or repelling maps, or linear maps with matrix has been considered. In this paper, we show Hyers-Ulam stability of difference equations induced by time dependent maps. Hyers-Ulam stability is abbreviated as HU stability. 
\medskip
Let $ F \colon \N \times \C \rightarrow \C $ be a function such that 
$$ |F_n(u)-F_n(v)| \leq p_n|u-v| $$
for every $ n \in \N $, $ u, v \in \C $. The positive numbers $ p_n $ is called {\bf contracting rate} of $ F_n $. Similarly, if the maps satisfy the inequalities
$$ |F_n(u)-F_n(v)| \geq p_n|u-v| , $$
$ p_n $ is called {\bf expanding rate}. Both contracting and expanding rate are called {\bf growth rate} of $ F_n $. The {\bf limit of the geometric average of growth rate} of all $ F_n $ is defined as
$$ \lim_{n \rightarrow \infty} \left( \prod_{j=1}^n p_j \right)^{\frac{1}{n}} $$
where all of $ p_n $ is either contracting or expanding rate of $ F_n $. Before taking the limit, the above product is called (geometric) {\bf average growth rate}. If the limit exists, then the limit is called {\bf exponential growth rate}. 
\noindent Note that if the contracting rates of all $ F_n $ is bounded above by $ p $ less than one, then difference equation has Hyers-Ulam stability. 
\begin{lem}[Corollary3 in \cite{jungnam}] \label{lem-HU stability in inv region}
Let $ F \colon \N \times \C \rightarrow \C $ be a function which satisfies that 
$$ |F(n,u)-F(n,v)| \leq p|u-v| $$
for all $ n \in \N $, $ u, v \in \C $ and for $ 0<p<1 $. For a given $ \e>0 $, suppose that the complex valued sequence $ \{a_n\}_{n\in \N} $ satisfies the inequality
\begin{align*} 
|a_{n+1} -F(n,a_n)| \leq \e
\end{align*}
for all $ n \in \N $. Then the sequence $ \{b_n\}_{n\in\N} $ defined as $ b_{n+1} = F(n,b_n) $ satisfies that 
$$ |b_n -a_n| \leq {p^{n-1}}|b_1-a_1| + \frac{1-{p^{n-1}}}{1-{p}}\e $$
for $ n \in \N $. 
\end{lem}

%
%

Contracting or expanding maps on higher dimension has similar properties. 
For various shadowing properties of uniformly contracting or expanding time varying maps, see \cite{sar}. For exponential trichotomy of infinite dimensional maps, see \cite{BD}.
However, the condition
$$ |F(n,u)-F(n,v)| \leq p|u-v| $$ 
for $ 0 < p <1 $ in Lemma \ref{lem-HU stability in inv region} is very sufficient for Hyers-Ulam stability. In this paper, we show Hyers-Ulam stability with exponential growth rate in one dimension. If the limit of the geometric average is exponential, the difference equation has Hyers-Ulam stability even though each time dependent map may have various contracting (or expanding) rate. This paper is a continuation of \cite{nam2}.

\section{HU stability for exponential growth rate less than one} \label{section 2}

\begin{lem} \label{lem-key lemma}
Let $ p_n $ be positive numbers for all $ n \in \N $. Let $ F \colon \N \times \C \rightarrow \C $ be a function which satisfies that 
$$ |F(n,u)-F(n,v)| \leq p_n|u-v| $$
for $ u, v \in \C $ and for all $ n \in \N $. For a given $ \e>0 $, let $ \{a_n\}_{n\in \N} $ be a complex valued sequence which satisfies the inequality
\begin{align*}  
|a_{n+1} -F(n,a_n)| \leq \e
\end{align*}
for all $ n \in \N $. Then the sequence $ \{b_n\}_{n\in\N} $ defined as $ b_{n+1} = F(n,b_n) $ satisfies that 
$$ |b_n -a_n| \leq  \left( \prod_{j=1}^{n-1} p_j \right)|b_1-a_1| + \left( \sum_{j=1}^{n-1} \prod_{i=j+1}^{n-1} p_i \right) \e $$
for $ n \in \N $.
\end{lem}

\begin{proof}
Suppose for induction that 
$$ |b_{n-1}-a_{n-1}| \leq \left( \prod_{j=1}^{n-2} p_j \right)|b_1-a_1| + \left( \sum_{j=1}^{n-2} \prod_{i=j+1}^{n-2} p_i \right) \e. $$
For $ n=1 $, the inequality trivially holds. The assumption and the above inequality imply that 
\begin{align*}
|b_n-a_n| &\leq |b_n-F(n-1,a_{n-1})| + |a_n-F(n-1,a_{n-1})| \\[0.3em]
&\leq |F(n-1,b_{n-1})-F(n-1,a_{n-1})| + |a_n-F(n-1,a_{n-1})| \\[0.3em]
&= p_{n-1}|b_{n-1}-a_{n-1}| + \e \\[0.3em]
&\leq p_{n-1} \left\{ \left( \prod_{j=1}^{n-2} p_j \right)|b_1-a_1| + \left( \sum_{j=1}^{n-2} \prod_{i=j+1}^{n-2} p_i \right) \e \right\} + \e \\[0.3em]
&= \left( \prod_{j=1}^{n-1} p_j \right)|b_1-a_1| + \left( \sum_{j=1}^{n-2} \prod_{i=j+1}^{n-1} p_i \right) \e + \e \\[0.3em]
&= \left( \prod_{j=1}^{n-1} p_j \right)|b_1-a_1| + \left( \sum_{j=1}^{n-1} \prod_{i=j+1}^{n-1} p_i \right) \e .
\end{align*}
The proof is complete. 
\end{proof}

\medskip

\begin{lem}[Lemma 5.4 in \cite{nam2}] \label{lem-key lemma2}
Let $ \{t_n\}_{n \in \N} $ be the sequence of positive real numbers. Then the limit $ \displaystyle \lim_{n \rightarrow \infty} t_n^{\frac{1}{n}} = 1 $ is satisfied if and only if  
$$ \lim_{n \rightarrow \infty} \frac{t_nK^n}{\sum_{j=1}^{n-1}t_jK^j} = K-1 $$
for $ K > 1 $. 
\end{lem}
\medskip

\begin{dfn}
Let $ \{w_n\}_{n \in \N} $ be a sequence which satisfies the inequality
$$ |w_{n+1} - g_n(w_n)| \leq \e $$
for a given $ \e > 0 $ and for all $ n \in \N $. 
 The difference equation $ z_{n+1} = g_n(z_n) $ is Hyers-Ulam stable if there exists a sequence $ \{z_n\}_{n \in \N} $ as follows
\begin{align*}
z_{n+1} = g_n(z_n) \quad \text{and} \quad |w_n - z_n| \leq G(\e)
\end{align*}
for all $ n \in \N $ and $ G(\e) \rightarrow 0 $ as $ \e \rightarrow 0 $. 
\end{dfn}
\medskip

\begin{thm} \label{thm-suffi condition of HU stability}
Let $ p_n $ be positive numbers for all $ n \in \N $. Let $ F(n,z) $ be a function which satisfies that 
$$ |F(n,u)-F(n,v)| \leq p_n|u-v| $$
for $ u, v \in \C $ and for every $ n \in \N $. 
Suppose that 
$$ \displaystyle \lim_{n \rightarrow \infty} \left( \prod_{j=1}^n p_j \right)^{\frac{1}{n}} = \frac{1}{K} <1 . $$
Then the difference equation $ b_{n+1} = F(n,b_n) $ has Hyers-Ulam stability.
\end{thm}

\begin{proof}
The assumption is equivalent that
$$ \lim_{n \rightarrow \infty} \left( \prod_{j=1}^n p_j^{-1} \right)^{\frac{1}{n}} = K > 1 . 
$$
Thus we may assume that $ \prod_{j=1}^n p_j^{-1} = t_nK^n $ where $ \lim_{n \rightarrow \infty} t_n^{\frac{1}{n}}=1 . $
Lemma \ref{lem-key lemma2} implies that 
\begin{align*} 
& \lim_{n \rightarrow \infty}\; \cfrac{\prod_{j=1}^{n} p_j^{-1}}{\sum\limits_{j=1}^{n-1} \prod_{i=1}^{j} p_i^{-1}} = K - 1 > 0 .
\end{align*}
This equation is satisfied if and only if
\begin{align*} 
& \lim_{n \rightarrow \infty} \left( \prod_{j=1}^{n} p_j \right) \sum_{j=1}^{n-1} \prod_{i=1}^{j} p_i^{-1} = \lim_{n \rightarrow \infty} \sum_{j=1}^{n-1} \prod_{i=j+1}^{n} p_i = 
\frac{1}{K-1} > 0 .
\end{align*}
Recall that $ \prod_{i=j+1}^{n} p_i $ is a positive number for $ n \in \N $, that is, $ \sum_{j=1}^{n} \prod_{i=j+1}^{n} p_i $ is an increasing sequence. Thus we have the inequality as follows
\begin{align} \label{eq-inequality for pn}
\sum_{j=1}^{n} \prod_{i=j+1}^{n} p_i = \sum_{j=1}^{n-1} \prod_{i=j+1}^{n} p_i + 1 \leq \frac{1}{K-1} + 1 = \frac{K}{K-1} 
\end{align}
for all $ n \in \N $. Let $ \{a_n\}_{n\in \N} $ be a sequence which satisfies the inequality
\begin{align}  \label{eq-seq a-n}
|a_{n+1} -F(n,a_n)| \leq \e
\end{align}
for all $ n \in \N $. Lemma \ref{lem-key lemma} and the inequality \eqref{eq-inequality for pn} imply that 
\begin{align*}
|b_n -a_n| &\leq  \left( \prod_{j=1}^{n-1} p_j \right)|b_1-a_1| + \left( \sum_{j=1}^{n-1} \prod_{i=j+1}^{n-1} p_i \right) \e \\
&\leq  \left( \prod_{j=1}^{n-1} p_j \right)|b_1-a_1| +  \frac{K\e}{K-1} .
\end{align*}
For each sequence $ \{a_n\}_{n\in \N} $ which satisfies \eqref{eq-seq a-n}, define $ b_1 $ is the same as $ a_1 $. Hence, the difference equation $ b_{n+1} = F(n,b_n) $ has Hyers-Ulam stability.  
\end{proof}

\bigskip

\section{HU stability for exponential growth rate greater than one }

\begin{dfn} \label{dfn-q-n}
Let $ q_n $ be complex numbers for $ n \in \N $. Let $ F \colon \N \times \C \rightarrow \C $ be a function. Define $ q_n $ as follows  
$$ \frac{F(n,u)-F(n,v)}{u-v} = q_n(u,v) $$
for $ u, v \in \C $ where $ u \neq v $. If $ F(n,\cdot\,) $ is differentiable at $ u \in \C $, then $ q_n(u,u) $ is defined as follows
$$ q_n(u,u) = \lim_{v \rightarrow u} \frac{F(n,u)-F(n,v)}{u-v} $$
for $  u \in \C $ and for $ n \in \N $. 
\end{dfn}

\noindent Observe that $ q_n(u,v) = q_n(v,u) $ for all $ u, v \in \C $. If $ F(n,\cdot\,) $ is differentiable on the set $ S \subset \C $, then $ q_n $ is continuous on $ S \times S $.  

\medskip

\begin{lem} \label{lem-key lemma 2nd}
Let $ \{a_n\}_{n\in \N} $ be a complex valued sequence which satisfies the equation
\begin{align*}  
a_{n+1} -F(n,a_n) = r_{n}
\end{align*}
%
for every $ n \in \N $. Then the sequence $ \{b_n\}_{n\in\N} $ defined by the equation $ b_{n+1} = F(n,b_n) $ for every $ n \in \N $ satisfies that 
$$ b_n -a_n =  \left( \prod_{j=1}^{n-1} q_j \right) (b_1-a_1) - \sum_{j=1}^{n-1}r_j \prod_{i=j+1}^{n-1} q_i $$
where $ q_j = q_j(b_j,a_j) $ for each $ j \in \N $.
\end{lem}

\begin{proof}
For $ n=1 $, the equation trivially holds. Suppose for induction that 
$$ b_{n-1}-a_{n-1} = \left( \prod_{j=1}^{n-2} q_j \right)(b_1-a_1) - \sum_{j=1}^{n-2}r_j \prod_{i=j+1}^{n-2} q_i . $$
The assumption of the sequence $ \{a_n\}_{n\in \N} $ and above equation imply that 
\begin{align*}
b_n-a_n &= b_n-F(n-1,a_{n-1}) - a_n + F(n-1,a_{n-1}) \\[0.3em]
&= F(n-1,b_{n-1})-F(n-1,a_{n-1}) - a_n + F(n-1,a_{n-1}) \\[0.3em]
&= q_{n-1}(b_{n-1}-a_{n-1}) - r_{n-1} \\[0.3em]
&= q_{n-1} \left\{ \left( \prod_{j=1}^{n-2} q_j \right)(b_1-a_1) - \sum_{j=1}^{n-2}r_j \prod_{i=j+1}^{n-2} q_i \right\} - r_{n-1} \\[0.3em]
&= \left( \prod_{j=1}^{n-1} q_j \right)(b_1-a_1) - \sum_{j=1}^{n-2}r_j \prod_{i=j+1}^{n-1} q_i - r_{n-1} \\[0.3em]
&= \left( \prod_{j=1}^{n-1} q_j \right)(b_1-a_1) - \sum_{j=1}^{n-1}r_j \prod_{i=j+1}^{n-1} q_i .
\end{align*}
The proof is complete. 
\end{proof}

\medskip


\noindent Denote $ F_{j} \circ F_{j-1} \circ \cdots \circ F_1 $ by $ {\bf F}_{j} $ where $ F(j,z) $ is expressed as $ F_j(z) $ for $ j \in \N $. Define $ {\bf F}_{0} $  as the identity map, that is, $ {\bf F}_{0} = \mathrm{Id} $. 
\medskip


\begin{lem} \label{lem-convergence and fixed point}
Let $ F_n $ be functions on $ \C $ 
for every $ n \in \N $. For a given $ \e > 0 $, let $ \{a_n\}_{n \in \N} $ is a given sequence 
defined as follows
\begin{align*}  
a_{n+1} -F_n(a_n) = r_{n}
\end{align*}
where $ |r_n| \leq \e $ for all $ n \in \N $. Let $ \{z_n\}_{n \in \N} $ is a given sequence defined by
$$ z_{n+1} = F_n(z_n) $$ 
for $ n \in \N $. Assume that there exists $ p_n > 0 $ which satisfies $ |F_n(u) -F_n(v)| \geq p_n|u -v| $ for $ u, v \in \C $ where the limit as follows
$$ \displaystyle \lim_{n \rightarrow \infty} \left( \prod_{j=1}^n  p_j \right)^{\frac{1}{n}} = {K} > 1 . $$
Then the following sum
\begin{align} \label{eq-seq convergent}
\sum_{j=1}^{\infty} \prod_{i=1}^j \frac{|a_i - z_i |}{| F_i(a_i) - F_i(z_i)|}
\end{align}
is convergent. The function 
$$ Q(z) = a_1 + \sum_{j=1}^{\infty} r_j\prod_{i=1}^j \frac{a_i - {\bf F}_{i}(z)}{F_i(a_i) -  {\bf F}_{i}(z)} $$
has the fixed point for all sufficiently small $ \e > 0 $.
\end{lem}
\begin{proof}
Firstly, we show the convergence of the series \eqref{eq-seq convergent}. The assumption yields the inequality
\begin{align} \label{eq-inequality for HU stability}
\frac{|a_i - z_i |}{| F_i(a_i) - F_i(z_i)|} \leq \frac{1}{p_i} .
\end{align}
Observe that 
\begin{align*} 
\lim_{n \rightarrow \infty} \left( \prod_{j=1}^n  \frac{1}{p_j} \right)^{\frac{1}{n}} = \frac{1}{K} < 1 .
\end{align*}
Thus comparison and root tests implies the series \eqref{eq-seq convergent} is convergent. \\
Secondly, the existence of the fixed point of $ Q(z) $ results from the contraction mapping theorem. Denote $ \frac{a_i - {\bf F}_{i}(z)}{F_i(a_i) -  {\bf F}_{i}(z)} $ by $ A_i $ temporarily. Define the norm of each $ A_i $ as 
$$ \| A_i \| = \sup_{z\neq a_i} \frac{| a_i - {\bf F}_{i}(z)|}{|F_i(a_i) -  {\bf F}_{i}(z)|} . $$
Since the inequality \eqref{eq-inequality for HU stability} is independent of the point $ z_i $, the inequality $ \| A_i \| \leq \frac{1}{p_i} $ holds for every $ i \in \N $. Thus we obtain that 
\begin{align*}
| Q(u) - Q(v) | &= \left|\,\sum_{j=1}^{\infty}r_j \prod_{i=1}^j A_i(u) - \sum_{j=1}^{\infty}r_j \prod_{i=1}^j A_i(v) \,\right| \\
&\leq \e \sum_{j=1}^{\infty}\prod_{i=1}^j \| A_i \| \cdot | u-v | \\
&\leq \e \sum_{j=1}^{\infty}\prod_{i=1}^j \frac{1}{p_i} \,| u-v | .
\end{align*}
The series $ \sum_{j=1}^{\infty}\prod_{i=1}^j \frac{1}{p_i} $ is convergent by root test. Then for sufficiently small $ \e > 0 $, the function $ Q $ is a contracting map. Hence, the contraction mapping theorem implies the existence of the unique fixed point of $ Q $. 
\end{proof}
\bigskip
\noindent HU stability depends on the maps $ {F_n} $ and the given $ \e > 0 $. In many cases $ \e $ must be small enough positive number. In addition to these conditions, HU stability may depend also on the set which contains the sequence $ \{a_n\} $ or $ a_1 $. For example, if the sequence of the maps $ F_n(z) = z^{n_k} $ where $ n_k $ is the positive integers for all $ k \in \N $ and $ \sup_k \{n_k\} < \infty $. Assume that $ \lim_{n \rightarrow \infty}\prod_{k=1}^n n_k = K>1 $.  HU stability is expected due to Theorem \ref{thm-suffi condition of HU stability 2}. However, HU stability is satisfied only on the exterior of the disk $ \{ z \colon |z| \geq \rho \} $ for a big enough $ \rho > 0 $. \\  
Let $ S $ be a set which contains the sequence $ \{a_n\}_{n \in \N} $ for all $ r_j $ and for a given $ \e > 0 $. Assume also that $ \displaystyle \inf_{j \in \N} \rm{dist}(a_j, \partial S) > 0 $.

\begin{thm} \label{thm-suffi condition of HU stability 2}
Let $ F_n $ be a complex valued function 
for each $ n \in \N $. For a given $ \e > 0 $, let $ \{a_n\}_{n \in \N} $ is a given sequence 
defined as follows
\begin{align*}  
a_{n+1} -F_n(a_n) = r_{n}
\end{align*}
where $ |r_n| \leq \e $ for all $ n \in \N $. Let $ q_n(u,v) $ be complex numbers for every $ n \in \N $ defined in Definition \ref{dfn-q-n}. 
Suppose that $ \displaystyle p_j \leq  \inf_{z \in S}  |q_j(a_j,z)| $ and the following limit
$$ \lim_{n \rightarrow \infty} \left( \prod_{j=1}^n  p_j \right)^{\frac{1}{n}} = {K} > 1 $$
is satisfied. Then the difference equation $ b_{n+1} = F(n,b_n) $ has Hyers-Ulam stability.
\end{thm}

\begin{proof} 
Lemma \ref{lem-key lemma 2nd} implies that 
$$ b_n -a_n =  \left( \prod_{j=1}^{n-1} q_j \right) (b_1-a_1) - \sum_{j=1}^{n-1}r_j \prod_{i=j+1}^{n-1} q_i $$
where $ q_n = q_n(b_n,a_n) $ for $ n \in \N $.
Then Lemma \ref{lem-convergence and fixed point} implies that $ b_1 $ can be chosen as follows 
\begin{align*} 
b_1 = a_1 + \sum_{j=1}^{\infty}r_j\prod_{i=1}^j \frac{1}{q_i}
\end{align*}
where $ q_i = {q_i(b_i,a_i)} $. The difference $ b_n -a_n $ is 
as follows 
\begin{align}
b_n -a_n &= \left(\prod_{j=1}^{n-1}q_j\right) \sum_{j=1}^{\infty} r_j \prod_{i=1}^{j}\frac{1}{q_i} - \sum_{j=1}^{n-1} r_j \prod_{i=j+1}^{n-1} q_i  \nonumber \\[0.5em]
&= \left( \prod_{j=1}^{n-1}q_j \right) \left( \sum_{j=1}^{\infty} r_j \prod_{i=1}^{j} \frac{1}{q_i} - \sum_{j=1}^{n-1} r_j \prod_{i=1}^{j} \frac{1}{q_i} \right) \nonumber \\[0.5em]
&= \left( \prod_{j=1}^{n-1}q_j \right) \left( \sum_{j=n}^{\infty} r_j \prod_{i=1}^{j}\frac{1}{q_i} \right) \nonumber \\[0.5em]
&= \sum_{j=n}^{\infty} r_j \prod_{i=n}^{j} \frac{1}{q_i} . \nonumber
\end{align}
Then the absolute value $ |b_n -a_n| $ is bounded as follows
\begin{align}
|b_n -a_n| &\leq \e \sum_{j=n}^{\infty} \left|\, \prod_{i=n}^{j}\frac{1}{q_i} \,\right| \nonumber \\[0.5em]
&\leq \e \sum_{j=n}^{\infty} \left|\, \prod_{i=n}^{j} \frac{1}{p_{i}} \,\right| = \e  \sum_{j=n}^{\infty} \left|\,\prod_{j=1}^{n-1} p_j \prod_{i=1}^{j} \frac{1}{p_{i}} \right| \nonumber \\[0.5em]
&= \e \left|\,\prod_{j=1}^{n-1} p_j \,\right| \left(\sum_{j=n}^{\infty} \left|\, \prod_{i=1}^{j} \frac{1}{p_{i}} \,\right| \right) \label{eq-bounds of diff of sequence2} \\
& = \e \left(t_{n-1}K^{n-1} \right) \sum_{j=n}^{\infty} \frac{1}{t_j K^j} .\label{eq-bounds of diff of sequences}
\end{align}
Let $ L $ is a positive number in $ (1,K) $, that is, $ K > L > 1 $. Denote $ t_n\left( \frac{L}{K} \right)^{n} $ by $ s_n $ for $ n \in \N $. Since the inequality $ \dfrac{1}{K^{j-n+1}} <  \dfrac{1}{L^{j-n+1}} $ is satisfied for all $ j \geq n-1 $, the following inequality holds
\begin{align*}
\sum_{j=n}^{\infty}\frac{t_{n-1}K^{n-1}}{t_j K^j} &< \sum_{j=n}^{\infty}\frac{t_{n-1}L^{n-1}}{t_j L^j} \\
&= t_{n-1}\left( \frac{L}{K} \right)^{n-1}K^{n-1} \sum_{j=n}^{\infty} \frac{1}{t_j \left( \frac{L}{K} \right)^{j}K^j} \\
&= s_{n-1}K^{n-1} \sum_{j=n}^{\infty} \frac{1}{s_j K^j} .
\end{align*}
Observe that $ \lim_{n \rightarrow \infty} s_n  = 0 $ and $ \lim_{n \rightarrow \infty} (s_n)^{\frac{1}{n}}  = \frac{L}{K} < 1 $. 
In order to estimate $ \e \left(t_{n-1}K^{n-1} \right) \sum_{j=n}^{\infty} \frac{1}{t_j K^j} $, we define $ C^1 $ function on the set of positive real numbers, say $ g $, which satisfies the following properties.
\begin{enumerate}
\item $ g(x) > 0 $ for all $ x > 0 $,
\item $ g(j) = \frac{1}{s_j} $ for every $ j \in \N $,
\item $ \lim_{n \rightarrow \infty} g(x) = \infty $, and 
\item $ \lim_{x \rightarrow \infty} \{g(x)\}^{\frac{1}{x}} = \frac{K}{L} > 1 $.
\end{enumerate}
The forth property is equivalent to the limit $ \lim_{x \rightarrow \infty} \frac{\ln(g(x))}{x} = \ln \left(\frac{K}{L} \right) $ or 
\begin{align} \label{eq-limit of g}
\lim_{x \rightarrow \infty} \frac{g'(x)}{g(x)} = \ln {K} - \ln {L} > 0
\end{align}
by L'Hospital's rule. Consider the function $ g(x)K^{-x} $ and calculate its derivative 
$$ \left( g(x)K^{-x}\right)' = \left[ \frac{g'(x)}{g(x)} - \ln K \right] g(x)K^{-x} . $$
The equation \eqref{eq-limit of g} implies that there exists big enough $ x_0 $ such that the inequality $ \frac{|g'(x)|}{g(x)} - \ln K < 0 $ for all $ x > x_0 $. Thus $ g $ is decreasing on the interval $ \{ x\ | \ x > x_0 \} $. Moreover, there exists a positive integer $ n $ as follows
\begin{align*}
\frac{|g'(x)|}{g(x)} \leq \frac{1}{2}\ln K 
\end{align*}
for all $ x \geq n-1 \geq x_0 $ with a suitable choice of $ L > \sqrt{K} $. 
Then we may use the comparison between infinite series and definite integral as follows 
\begin{align} \label{eq-comparison between series and integral}
\sum_{j=n}^{\infty} \frac{1}{s_j K^j} &\leq \int_{n-1}^{\infty}g(x)K^{-x}dx .
\end{align}
The definite integral on $ [n-1,\infty ) $ is estimated as follows 
\begin{align} \label{eq-bounds of integral}
\int_{n-1}^{\infty}g(x)K^{-x}dx 
&= \int_{n-1}^{\infty}g(x)e^{-x\ln K}dx \nonumber \\[0.5em]
&= \left[ -\frac{1}{\ln K}\, g(x)e^{-x\ln K} \right]_{n-1}^{\infty} + \frac{1}{\ln K} \int_{n-1}^{\infty}g'(x)e^{-x\ln K}dx \nonumber \\[0.5em]
&\leq \frac{1}{\ln K}\, g(n-1)e^{-(n-1)\ln K} + \frac{1}{2}\int_{n-1}^{\infty}g(x)K^{-x}dx .
\end{align}
Thus an upper bound of the integral \eqref{eq-bounds of integral} is as follows
\begin{align*}
\int_{n-1}^{\infty}g(x)K^{-x}dx \leq \frac{2}{\ln K}\,\frac{1}{s_{n-1}}K^{-(n-1)} .
\end{align*}
The estimations \eqref{eq-comparison between series and integral}  are applied to the difference $ |b_{n} -a_{n}| $ in the equation \eqref{eq-bounds of diff of sequences}. Then we have that
\begin{align*}
|b_{n} -a_{n}| &\leq \e \left(t_{n-1}K^{n-1} \right) \sum_{j=n}^{\infty} \frac{1}{t_j K^j} \\[0.2em]
&< \e \left(s_{n-1}K^{n-1} \right) \sum_{j=n}^{\infty} \frac{1}{s_j K^j} \\[0.2em]
&\leq \e \left(s_{n-1}K^{n-1} \right) \int_{n-1}^{\infty}g(x)K^{-x}dx \\[0.5em]
&= \e \left(s_{n-1}K^{n-1} \right) \frac{2}{\ln K}\,\frac{1}{s_{n-1}}K^{-(n-1)} \\
&= \frac{2\e}{\ln K} ,
\end{align*}
that is, we obtain an upper bound of difference between $ b_{n} $ and $ a_{n} $ as follows 
\begin{align} \label{eq-upper bound for HU stability}
|b_{n} -a_{n}| \leq \frac{2\e}{\ln K}
\end{align}
for all sufficiently large $ n \in \N $. An upper bound of $ |b_N -a_N | $ is estimated for $ N=1,2,\ldots,n-1 $ by the equation \eqref{eq-bounds of diff of sequence2} and the estimation \eqref{eq-upper bound for HU stability} as follows
\begin{align*}
|b_N-a_N| &\leq \e \sum_{j=N}^{\infty} \prod_{i=N}^{j} |q_{i}^{-1}| \leq \e \sum_{j=N}^{\infty} \prod_{i=N}^{j} |p_{i}^{-1}| \\
&= \e \left(\;\sum_{j=N}^{n} \prod_{i=N}^{j} |p_{i}^{-1}| + \sum_{j=n+1}^{\infty} \prod_{i=N}^{j} |p_{i}^{-1}| \right) \\[0.3em]
&= \e \left[\; \sum_{j=N}^{n} \prod_{i=N}^{j} |p_{i}^{-1}| + \left(\prod_{i=N}^{n} |p_{i}^{-1}| \right)\sum_{j=n+1}^{\infty} \prod_{i=n+1}^{j} |p_{i}^{-1}| \; \right] \\[0.3em]
&\leq \e \left[\; \sum_{j=N}^{n} \prod_{i=N}^{j} |p_{i}^{-1}| + \left(\prod_{i=N}^{n} |p_{i}^{-1}| \right) \frac{2}{\ln K}\right]
\end{align*} 
for $ N =1,2,\ldots,n-1 $. Hence, the difference equation $ b_{n+1} = F(n,b_n) $ has Hyers-Ulam stability.  
\end{proof}
\medskip

\begin{rk}
The unpper bound of $ | b_n - a_n | $ in Theorem \ref{thm-suffi condition of HU stability 2} would be better with a particular assumption of $ t_n $. For example, if we assume the boundedness of $ t_n $, that is, $ 0 < m = \inf_{n \in \N} \{ t_n \} \leq \sup_{n \in \N} \{ t_n \} = M < \infty $, then an upper bound of $ |b_n - a_n| $ is as follows 
\begin{align*}
|b_{n} -a_{n}| &\leq \e \left(t_{n-1}K^{n-1} \right) \sum_{j=n}^{\infty} \frac{1}{t_j K^j} \\[0.2em]
&\leq \e M K^{n-1}\sum_{j=n}^{\infty} \frac{1}{m K^j} \\
&= \frac{M\e}{m}K^{n-1} \frac{1}{K^n}\cdot\frac{1}{1-\frac{1}{K}} \\
&= \frac{M\e}{m(K-1)} ,
\end{align*}
which is a smaller upper bound than that of \eqref{eq-upper bound for HU stability} for large $ K > 1 $. 
\end{rk}

\section{No HU stability for periodic average growth rate}
The sequence $ \{a_n\} $ is {\it preperiodic} where there is a positive integer $ p $ such that $ a_{l+m} = a_l $ for every $ l \geq N $. If $ N =1 $, then the sequence is called the periodic sequence. The number $ m $ is said to be the period of sequence.

\begin{lem} \label{lem-suffi condition of HU stability 2}
Let $ p_n $ be positive numbers for every $ n \in \N $. Let $ F(n,z) $ be a function which satisfies that 
$$ |F(n,u)-F(n,v)| \leq p_n|u-v| $$
for $ u, v \in \C $ and for every $ n \in \N $. 
Suppose that the values $ \left( \prod_{j=1}^n p_j \right)^{\frac{1}{n}} $ are preperiodic and the maximum of the values is less than one. Then the difference equation $ b_{n+1} = F(n,b_n) $ has no Hyers-Ulam stability.
\end{lem}

\begin{proof}
Let 
$ \frac{1}{K_1}, \frac{1}{K_2},\ldots, \frac{1}{K_m} $ be all possible values of $ \left( \prod_{j=1}^n p_j \right)^{\frac{1}{n}} $ for every $ n > N $. Let $ m \geq 2 $ be the period of the sequence 
$ \left\{ \frac{1}{K_n} \right\}_{n\geq 0} $. Denote $ n $ by $ km + l $ where $ l = 1,2,\ldots,m $ for some $ k \in \N $. Let $ k_1 $ be the least positive number which satisfies that $ k_1m + l > N $. Thus $ K_l = K_{km+l} $ for every $ k = 0,1,2,\ldots $. Then we may assume that 
$$ \prod_{j=N+1}^{km + l} p_j = 
\frac{1}{(K_l)^{km + l-N}} $$
for every non negative integer $ k $ satisfying $ km + l > N $. 
We claim that $ \sum_{j=1}^{n} \prod_{i=j+1}^{n} p_i $ is unbounded for $ n \in \N $. It suffice to consider $ \left( \prod_{j=1}^n p_j \right)^{\frac{1}{n}} $ is periodic, that is, $ N = 0 $, because the first finite terms of the series does not affect convergence or divergence of the given series. Denote 
$  \max
\left\{ \frac{1}{K_1}, \frac{1}{K_2},\ldots, \frac{1}{K_m} \right\} $ by $ \frac{1}{K_p} $ which is less than one. Thus there exists $ \frac{1}{K_q} $  less than $ \frac{1}{K_p} $ where $ p \neq q $. 
The sum, $ \sum_{j=1}^{n} \prod_{i=j+1}^{n} p_i $, say $ s_n $, where $ n = km + p $ for $ k \in \N $ is estimated by the rearrangement of each terms as follows
\begin{align} %
\sum_{j=1}^{km+p} \prod_{i=j+1}^{km+p} p_i &= \sum_{j=1}^{km+p} \frac{(K_j)^j}{(K_p)^{km+p}} = \frac{1}{{(K_p)^{km+p}}}\sum_{j=1}^{km+p}{(K_j)^j} \nonumber \\
&= \frac{1}{{(K_p)^{km+p}}} \left[\,\sum_{j=0}^k(K_1)^{jm+1} + \sum_{j=0}^k(K_2)^{jm+2} + \cdots + \sum_{j=0}^k(K_p)^{jm+p} \right. \nonumber \\ &\qquad \left. + \sum_{j=0}^{k-1}(K_{p+1})^{jm+p+1} + \cdots + \sum_{j=0}^{k-1}(K_m)^{jm+m} \right] \nonumber \\
&\geq \frac{1}{{(K_p)^{km+p}}}\sum_{j=0}^{k-1}{(K_q)^{jm+q}} \nonumber \\
&= \frac{1}{{(K_p)^{km+p}}}\cdot \frac{(K_q)^q((K_q)^{km}-1)}{(K_q)^m-1} \nonumber \\
&> \frac{1}{{(K_p)^{km+p}}}\cdot  \frac{(K_q)^q (K_q)^{km}}{(K_q)^m} \nonumber \\
&= \left(\frac{K_q}{K_p}\right)^{km} \frac{(K_q)^{q}}{(K_p)^{p+m}} . \label{eq-bound of prod of p}
\end{align}
The fact that $ K_q > K_p $ implies the above inequality \eqref{eq-bound of prod of p} is unbounded for $ k \in \N $. Since a subsequence $ \{s_{km+p}\}_{k\in \N} $ of $ \{s_n\}_{n\in \N} $ is unbounded, so is the sequence $ \{s_n\}_{n\in \N} $. Let $ \{a_n\}_{n\in \N} $ be a sequence which satisfies the inequality
\begin{align*}  
a_{n+1} -F(n,a_n) = r_n
\end{align*}
for every $ n \in \N $ where $ |r_n| \leq \e $ for all $ n \in \N $. Let $ \{b_n\}_{n\in\N} $ be the sequence defined by the equation $ b_{n+1} = F(n,b_n) $ for every $ n \in \N $. Lemma \ref{lem-key lemma 2nd} implies that  
$$ b_n -a_n =  \left( \prod_{j=1}^{n-1} q_j \right) (b_1-a_1) - \sum_{j=1}^{n-1}r_j \prod_{i=j+1}^{n-1} q_i $$
where $ q_j = q_j(b_j,a_j) $ for each $ j \in \N $. We may choose $ r_n = \e $ for all $ n \in \N $. Then

\begin{align} \label{eq-triangular ineq of difference}
|b_n -a_n| + \left| \prod_{j=1}^{n-1} q_j \right|\,|b_1-a_1| &\geq \left( \sum_{j=1}^{n-1} \prod_{i=j+1}^{n-1} p_i \right) \e .
\end{align}
Observe that  
$$ \left| \prod_{j=1}^{n-1} q_j \right| \leq \left| \prod_{j=1}^{n-1} p_j \right| \leq \left(\frac{1}{K_p}\right)^{n-1}  < 1 , $$
which is bounded for all $ n \in \N $. Thus for any choice of $ b_1 $, the second term of the left side of equation \eqref{eq-triangular ineq of difference} is bounded. However, the equation \eqref{eq-bound of prod of p} implies that $ \left| \sum_{j=1}^{n-1} \prod_{i=j+1}^{n-1} p_i \right| $ is unbounded for $ n \in \N $. Then $ |b_n -a_n| $ is also unbounded for $ n \in \N $. Hence, the difference equation $ b_{n+1} = F(n,b_n) $ has no Hyers-Ulam stability.  
\end{proof}

\medskip
%

\begin{thm}
\label{thm-suffi condition of HU stability 3}
Let $ p_n $ be positive numbers for all $ n \in \N $. Let $ F(n,z) $ be a function which satisfies that 
$$ |F(n,u)-F(n,v)| \leq p_n|u-v| $$
for $ u, v \in \C $ and for every $ n \in \N $. Denote $ n $ by $ km + l $ where $ k $ is a non-negative number, $ m \geq 2 $ is some positive integer and $ l = 1,2,\ldots,m $. 
Suppose that the values $ \left( C_n\prod_{j=1}^n p_j \right)^{\frac{1}{n}} $ for $ n \in \N $ is preperiodic with period $ m $ where positive constants $ C_{km+l} = C_l $ is for all non-negative integer $ k $. Suppose also that the maximum of preperiodic values is less than one. Then the difference equation $ b_{n+1} = F(n,b_n) $ has no Hyers-Ulam stability.
\end{thm}

\begin{proof}
The proof is similar to that of Lemma \ref{lem-suffi condition of HU stability 2}. We use the same notions for the proof. Thus we obtain that
\begin{align} \label{eq-product of periodic constant2}
C_l\prod_{j=N+1}^{km + l} p_j = \frac{1}{(K_l)^{km + l-N}}
\end{align}
for every $ k = 0,1,2,\ldots $. 
Let the numbers $ \frac{1}{K_1}, \frac{1}{K_2},\ldots,\frac{1}{K_m} $ be all possible positive values of \eqref{eq-product of periodic constant2} and denote the maximum of these numbers by $ \frac{1}{K_p} $. Since $ m \geq 2 $, there exists $ K_q > K_p $ for some $ q \neq p $. We claim that $ \sum_{j=N+1}^{km+p} \prod_{i=j+1}^{km+p} p_i $ is unbounded for $ k $. We may assume that the preperiodic values is periodic, that is, $ N=0 $ because the first finite terms of the partial sums does not affect the convergence or divergence of the given series. The similar calculation of \eqref{eq-bound of prod of p} in Lemma \ref{lem-suffi condition of HU stability 2} yields the following estimations
\begin{align} 
\sum_{j=1}^{km+p} \prod_{i=j+1}^{km+p} p_i &= \frac{1}{C_p} \sum_{j=1}^{km+p} \frac{(K_j)^j}{(K_p)^{km+p}} = \frac{1}{C_p}\frac{1}{{(K_p)^{km+p}}}\sum_{j=1}^{km+p}{(K_j)^j} \nonumber \\
&\geq \frac{1}{{C_p(K_p)^{km+p}}}\sum_{j=0}^{k-1}{(K_q)^{jm+q}} \nonumber \\
&> \frac{1}{{C_p(K_p)^{km+p}}}\cdot  \frac{(K_q)^q (K_q)^{km}}{(K_q)^m} \nonumber \\
&= \left(\frac{K_q}{K_p}\right)^{km} \frac{(K_q)^{q}}{C_p(K_p)^{p+m}} , \label{eq-bound of prod of p1}
\end{align}
which is unbounded for $ k \in \N $. For any given $ \e > 0 $, let $ \{a_n\}_{n\in \N} $ be a sequence which satisfies the inequality
\begin{align*}  
a_{n+1} -F(n,a_n) = \e
\end{align*}
for all $ n \in \N $. Lemma \ref{lem-key lemma 2nd} and triangular inequality imply that 
\begin{align} \label{eq-triangular inequality}
|b_n -a_n| + \left| \prod_{j=1}^{n-1} q_j \right|\,|b_1-a_1| \geq  \left( \sum_{j=1}^{n-1} \prod_{i=j+1}^{n-1} p_i \right) \e 
\end{align}
where $ q_j = q_j(a_j,b_j) $ for each $ j \in \N $. The right side of the equation \eqref{eq-triangular inequality} is unbounded for $ n \in \N $ by the estimation \eqref{eq-bound of prod of p1}. 
Hence, the difference equation $ b_{n+1} = F(n,b_n) $ has no  Hyers-Ulam stability.
\end{proof}

\section{Applications}

\begin{exa} \label{exa-exponential growth less than one}
Let $ f(n,z) $ be the functions as follows
\begin{align*}
f(n,z) = 
\begin{cases}
2z, & n = 1,3,5,\ldots \\
\frac{1}{3}z, & n= 2,4,6, \ldots
\end{cases} .
\end{align*}
Thus $ |f(n,u)-f(n,v)| \leq p_n |u-v| $ where 
\begin{align*}
p_n =  
\begin{cases}
2, & n = 1,3,5,\ldots \\
\frac{1}{3}, & n= 2,4,6, \ldots
\end{cases} .
\end{align*}
Thus we obtain that $ \displaystyle \lim_{n \rightarrow \infty} \left( \prod_{j=1}^n p_j \right)^{\frac{1}{n}} = \sqrt{\frac{2}{3}} $ because 
$$ \lim_{k \rightarrow \infty} \left( \prod_{j=1}^{2k} p_j \right)^{\frac{1}{2k}} = \lim_{k \rightarrow \infty} \left( \prod_{j=1}^{k} \frac{2}{3} \right)^{\frac{1}{2k}} = \sqrt{\frac{2}{3}} $$
for $ n = 2k,\ k=1,2,3,\ldots $ and moreover
$$ \lim_{k \rightarrow \infty} \left( \prod_{j=1}^{2k+1} p_j \right)^{\frac{1}{2k+1}} = \lim_{k \rightarrow \infty} \left( 2 \prod_{j=1}^{k} \frac{2}{3} \right)^{\frac{1}{2k+1}} = \sqrt{\frac{2}{3}}
$$
for $ n = 2k+1,\ k=0,1,2,\ldots $. Theorem \ref{thm-suffi condition of HU stability} implies that the difference equation $ b_{n+1} = f(n,b_n) $ has Hyers-Ulam stability. Let us choose that $ b_1= a_1 $. Then we obtain that 
$$ | b_n - a_n | \leq \frac{\sqrt{\frac{3}{2}}\e}{\sqrt{\frac{3}{2}}-1} = \frac{\sqrt{3}\,\e}{\sqrt{3}-\sqrt{2}} = \sqrt{3}\left(\sqrt{3} + \sqrt{2} \right)\e =  \left(3+\sqrt{6} \right)\e <  6\e  $$
for all $ n \in \N $. Since $ |b_n - a_n| < 6\e $ for all $ n \in \N $, the difference equation $ b_{n+1} = f(n,b_n) $ has Hyers-Ulam stability. 
\end{exa}

\medskip

\begin{exa} \label{exa-exponential growth greater than one}
Let $ g(n,z) $ be the functions as follows
\begin{align*}
g(n,z) = 
\begin{cases}
3nz, & n = 1,3,5,\ldots \\
\frac{1}{2n}z, & n= 2,4,6, \ldots
\end{cases} .
\end{align*}
Thus $ |g(n,u)-g(n,v)| \leq p_n |u-v| $ where 
\begin{align*}
p_n =  
\begin{cases}
3n, & n = 1,3,5,\ldots \\
\frac{1}{2n}, & n= 2,4,6, \ldots
\end{cases} .
\end{align*}
Thus we obtain that 
\begin{align*}
\prod_{j=1}^n p_j = 
\begin{cases}
\dfrac{3^k}{2^k}\cdot \dfrac{1\cdot 3\cdot 5\cdots (2k-1)}{2\cdot 4\cdot 6\cdots (2k)}, & n = 2k,\qquad k= 1,2,3,\ldots \vspace{2mm} \\
\dfrac{3^{k+1}}{2^k}\cdot \dfrac{1\cdot 3\cdot 5\cdots (2k-1)(2k+1)}{2\cdot 4\cdot 6\cdots (2k)}, & n= 2k+1,\ k = 1,2,3, \ldots
\end{cases} .
\end{align*}
In order to calculate the limit of $ \prod_{j=1}^n p_j $, we use the following inequalities
\begin{align} \label{eq-product ineq1}
\frac{1}{\sqrt{4k+1}} &\leq \dfrac{1\cdot 3\cdot 5\cdots (2k-1)}{2\cdot 4\cdot 6\cdots (2k)} \leq \frac{1}{\sqrt{3k+1}} \\[0.5em]
\frac{3}{\sqrt{2}}\sqrt{2k+1} \leq \frac{3(2k+1)}{\sqrt{4k+1}} &\leq 3\,\dfrac{1\cdot 3\cdot 5\cdots (2k-1)}{2\cdot 4\cdot 6\cdots (2k)}(2k+1) \leq \frac{3(2k+1)}{\sqrt{3k+1}} \leq 3\sqrt{2k+1} \nonumber
\end{align}
for every $ k \in \N $. Observe that for any rational expression $ Q(k) $, the limit $ \lim_{k \rightarrow \infty} \left(\sqrt{Q(k)}\right)^{\frac{1}{k}} = 1 $ is satisfied. Since the squeezing of the limits as follows 
\begin{align*}
\lim_{k \rightarrow \infty} \left( \frac{1}{\sqrt{4k+1}}\left(\frac{3}{2} \right)^k \right)^{\frac{1}{2k}} & \leq \lim_{k \rightarrow \infty} \left( \prod_{j=1}^{2k} p_j \right)^{\frac{1}{2k}} \leq \lim_{k \rightarrow \infty} \left( \frac{1}{\sqrt{3k+1}}\left(\frac{3}{2} \right)^k \right)^{\frac{1}{2k}} \\[0.5em]
\lim_{k \rightarrow \infty} \left( \frac{3}{\sqrt{2}}\sqrt{2k+1} \left(\frac{3}{2} \right)^k  \right)^{\frac{1}{2k+1}} & \leq \lim_{k \rightarrow \infty} \left( \prod_{j=1}^{2k+1} p_j \right)^{\frac{1}{2k+1}} \leq \lim_{k \rightarrow \infty} \left(  3\sqrt{2k+1} \left(\frac{3}{2} \right)^k \right)^{\frac{1}{2k+1}} 
\end{align*}
we obtain the limit
$$ 
\lim_{n \rightarrow \infty} \left( \prod_{j=1}^n p_j \right)^{\frac{1}{n}} = \sqrt{\frac{3}{2}} . $$
Then Theorem \ref{thm-suffi condition of HU stability 2} implies that the difference equation $ b_{n+1} = g(n,b_n) $ has Hyers-Ulam stability. 

\end{exa}

\medskip

\begin{exa} 
Let $ h(n,z) $ be the functions as follows
\begin{align*}
h(n,z) = 
\begin{cases}
2^nz, & n = 1,3,5,\ldots \vspace{2mm} \\
\dfrac{1}{2^{n+3}}z, & n= 2,4,6, \ldots
\end{cases} .
\end{align*}
Thus $ |h(n,u)-h(n,v)| \leq p_n |u-v| $ where 
\begin{align*}
p_n =  
\begin{cases}
2^n, & n = 1,3,5,\ldots \vspace{2mm} \\
\dfrac{1}{2^{n+3}}, & n= 2,4,6, \ldots
\end{cases} .
\end{align*}
Thus we obtain the product of $ p_n  $ as follows by straightforward calculation  
\begin{align*}
\prod_{j=1}^n p_j = 
\begin{cases}
\dfrac{1}{2^{n-2}}\, & n = 1,3,5,\ldots \vspace{2mm} \\
\dfrac{1}{2^{2n}}, & n= 2,4,6, \ldots
\end{cases} .
\end{align*}
Then we have that 
\begin{align*}
\left( 4\prod_{j=1}^{2k-1} p_j \right)^{\frac{1}{2k-1}} = \frac{1}{2} \quad \text{and} \quad \left(\prod_{j=1}^{2k} p_j \right)^{\frac{1}{2k}} = \frac{1}{4}  
\end{align*}
for all $ k \in \N $. The maximum of the set $ \left\{ \left( 4 \prod_{j=1}^{2k-1}p_j \right)^{\frac{1}{2k-1}},\ \left(\prod_{j=1}^{2k} p_j \right)^{\frac{1}{2k}}  \right\} $ is $ \frac{1}{2} $, which is less than one. Hence, Theorem \ref{thm-suffi condition of HU stability 3} implies that the difference equation $ b_{n+1} = h(n,b_n) $ has no Hyers-Ulam stability. 
\end{exa}

\medskip

\begin{exa}
Let $ f(n,x) $ be the function from $ \R $ to $ \R $ defined as follows 
$$ f(n,x) = 3x + \frac{1}{n}\sin \left(\frac{x}{n}\right) $$
for $ n \in \N $. Thus $ q_n $ of $ f(n,x) $ defined in Definition \ref{dfn-q-n} is as follows
\begin{align*}
q_n(u,v) = 3 + \frac{1}{n^2}\frac{\sin \left(\frac{u}{n}\right) - \sin \left(\frac{v}{n}\right)}{\frac{u}{n} - \frac{v}{n}}
\end{align*}
for $ n \in \N $. Thus Mean Value Theorem implies that 
$$  q_n(u,v) \geq 3 - \frac{1}{n^2} 
$$ 
for every $ n \in \N $. Denote $ 3 - \frac{1}{n^2} $ by $ p_n $. Since the product of numbers $ \frac{1}{3^n}\prod_{j=1}^n \left(3 - \frac{1}{n^2}\right) $ is convergent as $ n \rightarrow \infty $, we obtain that 
$$ \lim_{n \rightarrow \infty} \left(\prod_{j=1}^n p_j \right)^{\frac{1}{n}} = 3 > 1 . $$ 
Hence, Theorem \ref{thm-suffi condition of HU stability 2} implies that the difference equation $ b_{n+1} = f(n,b_n) $ has Hyers-Ulam stability. 

\end{exa}

\section*{Note and discussion}
Equation \eqref{eq-product ineq1} in Example \ref{exa-exponential growth greater than one} can be shown by induction. The ratio of product of even and odd integers can be expressed as the ratio of two gamma functions. For various and refined evaluation, see \cite{QL}. \\
If the complex analytic function $ F $ satisfies $ |F(u) - F(v)| \leq C|u-v| $ for some $ C>0 $ and for all $ u, v \in \C $, then $ F $ is a polynomial of degree one by Liouville's theorem. Thus in most cases we require a domain $ S $ which is not $ \C $, which is (forward) invariant under all given functions, that is, $ F_n(S) \subset S $ for all $ n \in \N $. For example, let $ f $ and $ g $ be polynomials of which degree is greater than or equal to two. Assume that the filled Julia set of $ f $ contains that of $ g $. Thus we may choose the Fatou component of $ \infty $ as an invariant set under $ f $ and $ g $, say $ S $, rather than exterior of the disk with sufficiently big radius. If we choose $ F_n $ as the arbitrary composition of functions $ f $ and $ g $ for every $ n \in \N $, then the inequality $ |F_n(u) -F_n(v)| \geq p_n|u-v| $ for all $ n \in \N $ and $ \inf \{p_n\} > 1 $. The positive number $ \e $ may depend on the distance between $ a_1 $ and $ \partial S $. For the dynamics of polynomial semigroup, see \cite{SS}. \\
HU stability of difference equation induced one dimensional map may be generalized on higher dimension. For instance, the matrix difference equation is defined as follows
\begin{align*}
\begin{bmatrix}
z_{n+1} \\
w_{n+1}
\end{bmatrix} = 
\begin{bmatrix}
p_n & 0 \\
0 & q_n
\end{bmatrix}
\begin{bmatrix}
z_{n} \\
w_{n}
\end{bmatrix}
\end{align*}
where $ p_n $ and $ q_n $ are coefficients of maps defined in Example \ref{exa-exponential growth less than one} and Example \ref{exa-exponential growth greater than one}. Thus the above matrix difference equation has HU stability because difference equations $ z_{n+1} = f_n(z_n) $ and $ z_{n+1} = g_n(z_n) $ have HU stability. However, decomposed equations are not induced by either contracting or expanding maps. \\



\end{document}